\documentclass[a4paper,11pt]{article}
\pdfoutput = 0
\usepackage{stmaryrd}
\usepackage[normalem]{ulem}
\usepackage{graphicx}
\usepackage{theorem}
\usepackage{amsmath,amssymb,mathrsfs}
\usepackage[all]{xy}
\usepackage{color}
\usepackage{times}
\usepackage{enumerate,vmargin}
\usepackage{verbatim}
\usepackage{paralist} 
\usepackage[labelfont=bf,font=it]{caption}
\usepackage[numbers,sort&compress]{natbib}

\setpapersize{A4}
\setmarginsrb{2.5cm}{2cm}{2.5cm}{2cm}{.4cm}{4mm}{0cm}{7.5mm}

\newtheorem{lem}{Lemma}
\newtheorem{thm}[lem]{Theorem}

\newtheorem{cor}[lem]{Corollary}

{\theorembodyfont{\rmfamily}\newtheorem{rem}{Remark}}

\newcommand{\noi}{\noindent}

\newcommand{\llb}{\llbracket}
\newcommand{\rrb}{\rrbracket}

\newcommand{\cT}{\mathcal T}

\newcommand{\cF}{\mathcal F}

\newcommand{\cP}{\mathcal P}
\newcommand{\cQ}{\mathcal Q}
\newcommand{\cL}{\mathcal L}

\newcommand{\cI}{\mathcal I}
\newcommand{\cJ}{\mathcal J}

\newcommand{\bbR}{\mathbb R}
\newcommand{\bbZ}{\mathbb Z}

\newcommand{\bbN}{\mathbb N}

\newcommand{\bbC}{\mathbb C}

\newcommand{\bE}{\mathbf E}
\newcommand{\un}{\mathbf 1}
\newcommand{\bC}{\mathbf C}
\newcommand{\bN}{\mathbf N}

\newcommand{\bP}{\mathbf P}

\newcommand{\rT}{\mathrm T}

\renewcommand{\(}{\left(}
\renewcommand{\)}{\right)}

\newcommand{\lam}{\lambda}

\newcommand{\indi}{\boldsymbol{1}}

\def\cq{$\hfill \square$}
\def\cqfd{$\hfill \blacksquare$}

\title{{\huge \textsc{Height and Diameter of Brownian tree}}
\thanks{Research partly supported by ANR GRAAL Projet.}}
\author{Minmin \textsc{Wang}
\thanks{\textbf{Institution}: PRES Sorbonne Universit\'es, UPMC Universit\'e Paris 06,  LPMA (UMR 7599). \textbf{Postal address}: LPMA, Bo\^ite courrier 188, 4 place Jussieu, 75252 Paris Cedex 05, FRANCE. \textbf{Email}: wangminmin03@gmail.com} }

\date{\today}

\begin{document}

\maketitle


\begin{abstract}
By computations on generating functions, Szekeres proved in 1983 that the law of the diameter of a uniformly distributed rooted labelled tree 
with $n$ vertices, rescaled by a factor $n^{-\frac{1}{2}}$, converges to a distribution whose density is explicit. 
Aldous observed in 1991 that this limiting distribution is the law of the diameter 
of the Brownian tree. In our article, we provide a computation of this law which is directly based 
on the normalized Brownian excursion. Moreover, we 
provide an explicit formula for the joint law of the height and diameter of the Brownian tree, which is a new result. 
\smallskip


\smallskip

\end{abstract}


\section{Introduction}
For any integer $n\! \ge \! 1$, let $T_n$ be a uniformly distributed random rooted labelled tree with $n$ vertices and we denote by $D_n$ its diameter with respect to the graph distance. By computations on generating functions, \citet{Sz83} proved that 
\begin{equation}\label{eq: cv_dm}
n^{-\frac{1}{2}}D_n\overset{\textrm{(law)}}{-\!\!\! -\!\!\!\longrightarrow} \Delta  \; , 
\end{equation}
where $\Delta$ is a random variable whose probability density $f_\Delta$ is given by
\begin{equation}\label{eq: sz}
f_\Delta(y)=\frac{\sqrt{2\pi}}{3}\sum_{n\ge 1}\(\frac{64}{y^4}(4b^4_{n, y}-36b^3_{n, y}+75b^2_{n, y}-30b_{n, y})+\frac{16}{y^2}(2b^3_{n, y}-5b^2_{n, y})\) e^{-b_{n, y}} , 
\end{equation}
where $b_{n, y}\! :=\! 8(\pi n/y)^2$, for all $y \! \in \! (0, \infty)$ and for all integers $n\! \ge \! 1$. This result is  implicitly written in Szekeres \cite{Sz83} p.~395 formula (12). 
See also \citet{BrFl} for a similar result for binary trees. On the other hand, 
Aldous \cite{aldcrt1, aldcrt3} has proved that $T_n$, whose graph distance is rescaled by a factor $n^{-\frac{1}{2}}$, 
converges in distribution to the Brownian tree (also called Continuum Random Tree) that is a random compact metric space. From this, Aldous has deduced that $\Delta$ has the same distribution as the diameter of the Brownian tree: see \cite{aldcrt2}, Section 3.4, 
(though formula (41) there is not accurate). As proved by Aldous \cite{aldcrt3} and by Le Gall \cite{Legall1993}, the Brownian tree is coded by the normalized Brownian excursion of length 
$1$ (see below for more detail). Then, the question was raised by Aldous \cite{aldcrt2} that whether we can establish \eqref{eq: sz} directly from computations on the normalized Brownian excursion. (See also \cite{pitmanstflour}, Exercise 9.4.1.)
In this work, we present a solution to this question: we compute the Laplace transform for the 
law of the diameter of the Brownian tree based on Williams' decomposition of Brownian excursions. We also provide a formula for the joint law of the total height and diameter of the Brownian tree, which appears to be new. 
In a joint work with Duquesne \cite{DuWa14}, we generalize the present method to a L\'evy tree and study its total height and diameter. 
Before stating precisely our results, let us first recall the definition of the Brownian tree coded by the normalized Brownian excursion.

\paragraph{Normalized Brownian excursion.} Let $X\! =\! (X_t)_{t\ge 0}$ be a continuous process 
defined on a probability space $(\Omega,\cF, \bP)$ such that $(\frac{_1}{^{\sqrt{2}}}X_t)_{ t\ge 0}$ is distributed as a linear standard Brownian motion such that 
$\bP (X_0\! = \! 0)\! =1\! $ (the reason for the normalizing constant $\sqrt{2}$ is explained below). Thus, 
$$ \forall u \in \bbR, \; t\in \bbR_+, \quad \bE \big[ e^{iuX_t}\big]= e^{-tu^2} \; .$$
For all $t\! \in \! \bbR_+$, we set $I_t \! =\! \inf_{s\in [0, t]} X_s$. Then, the reflected process $X\! -\! I$ is a strong Markov process, the state $0$ is instantaneous in 
$(0,\infty)$ and recurrent, and $\!-I$ is a local time at level $0$ for $X\! -\! I$ (see Bertoin \cite{Bebook}, Chapter VI). 
We denote by $\bN$ the \textit{excursion measure} associated with the local time $\! -I$; $\bN$ is a sigma finite measure on the space of continuous paths $\bC (\bbR_+, \bbR_+)$. 
More precisely, let $\bigcup_{i\in \cI} (a_i, b_i) \! =\!  \big\{ t \! >\! 0: X_t \! -\! I_t  \! >\! 0 \big\}$ 
be the excursion intervals of the reflected process $X\!-\! I$ above $0$; 
for all $i\in \cI$, we set $e_i (s)\! = \! X_{(a_i +s)\wedge b_i} \! -\!  I_{a_i}$, $s\! \in \! \bbR_+$. 
Then, 
\begin{equation}
\label{Nexcdef}
\textrm{$\sum_{i\in \cI} \delta_{(-I_{a_i}, e_i)}$ is a Poisson point measure on $\bbR_+ \! \times \! \bC (\bbR_+, \bbR_+)$ with intensity $dt \, \bN(de)$.}
\end{equation} 
We shall denote by $e\! = \! (e_t)_{t\ge 0}$ the canonical process on $\bC(\bbR_+, \bbR_+)$. We define its \textit{lifetime} by 
\begin{equation}
\label{lifetimedef}
\zeta\! = \! \sup \{ t \! \ge \! 0: e_t \! >\! 0\} \; ,   
\end{equation}
with the convention that $\sup \varnothing \! = \! 0$. 
Then, $\bN$-a.e.~$e_0 \! = \! 0$, $\zeta \! \in \! (0,\infty)$ and for all $t\! \in \! (0,\zeta)$, $e_t\! >\! 0$. 
Moreover, one has 
\begin{equation}\label{loizeta}
\forall\lambda\in (0, \infty), \qquad \bN \big(1\! -\! e^{-\lambda \zeta} \big)= \sqrt{\lambda} \quad {\rm and}\quad
\bN \big(\zeta \! \in \! d r \big)= \frac{dr }{2\sqrt{\pi}\,  r^{3/2}} \; .
\end{equation}
See Blumenthal \cite{Blubook} IV.1 for more detail. 

Let us briefly recall the scaling property of $e$ under $\bN$.  
To that end, recall that $X$ satisfies the following scaling property: for all  $r \! \in \! (0, \infty)$, 
$(r^{-\frac{1}{2}}\! X_{rt} )_{ t\ge 0}$ 
has the same law as 
$X$, which easily entails that 
\begin{equation}
\label{HscalingN_br}
\big( r^{-\frac{1}{2}} e_{rt} \big)_{t\ge 0} \quad  \textrm{under} \quad  
r^{\frac{1}{2}} \,  \bN \quad  
\overset{\textrm{(law)}}{=} \quad e  \quad  \textrm{under}  \quad \bN \; . 
\end{equation}
This scaling property implies that there exists a family of laws on $\bC(\bbR_+, \bbR_+)$ denoted by  
$\bN(\, \cdot \, | \, \zeta \! = \! r)$, $r \! \in \! (0, \infty)$, such that 
$r \mapsto \bN(\, \cdot \, | \, \zeta \! =\!  r)$ is weakly continuous on $\bC(\bbR_+, \bbR_+)$, 
such that $\bN(\, \cdot \, | \, \zeta \! = \! r)$-a.s.~$\zeta \! = \! r$ and such that 
\begin{equation}
\label{zetades_br}
\bN = \int_0^\infty \!\!\! \bN(\, \cdot \, | \, \zeta \! =\!  r) \, \bN \big( \zeta \! \in \! d r \big) \; .
\end{equation}
Moreover, by (\ref{HscalingN_br}), $\big( r^{-\frac{1}{2}} e_{rt} \big)_{t\ge 0}$ under 
$\bN(\, \cdot \, | \, \zeta \! = \! r)$ has the same law as $e$ under 
$\bN(\, \cdot \,| \, \zeta \! =\!  1 )$. To simplify notation we set 
\begin{equation}
\label{Nnorm_br}
\bN_{\! {\rm nr}}:= \bN(\, \cdot \,| \, \zeta \! =\!  1 ) \; .
\end{equation}
Thus, for all measurable functions $F\! :\! \bC(\bbR_+, \bbR_+)\! \rightarrow \! \bbR_+$, 
\begin{equation}
\label{echtscal_br}
\bN \big[ F(e)\big]=  \frac{1}{2\sqrt{\pi}} \int_0^\infty \!\!\! \! dr  \, r^{-\frac{3}{2}}\, \bN_{\! \textrm{nr}}\Big[ F\Big( \big( r^{\frac{1}{2}} e_{t/r} \big)_{\! t\ge 0}\Big) \Big] \; . 
\end{equation}
\begin{rem}
\label{nonstandard} The standard Ito measure $\bN_{\! {\rm Ito}}^+$ of positive excursions, as defined for instance in Revuz \& Yor \cite{yorbook} Chapter XII Theorem 4.2, is derived from $\bN$ by the following scaling relations: 
$$ \textrm{$\bN_{\! {\rm Ito}}^+$ is the law of $\frac{_1}{^{\sqrt{2}}} e$ under $\frac{_1}{^{\sqrt{2}}}\bN$ and thus, $\bN_{\! {\rm Ito}}^+( \, \cdot \, | \, \zeta \! = \! 1)$ is the law of 
$\frac{_1}{^{\sqrt{2}}} e$ under $\bN_{\! {\rm nr}}$.} $$
Consequently, the law $\bN_{\! {\rm nr}}$ is \textit{not} the standard version for normalized Brownian excursion measure. However, we shall refer to it as the \textit{normalized Brownian excursion measure}.  \cq 
\end{rem}

\paragraph{Real trees.} Let us recall the definition of \textit{real trees} that are metric spaces generalizing graph-trees: 
let $(T, d)$ be a metric space; it is a real tree if the following statements hold true.

\smallskip

\noi
(a)~~For all $\sigma_1, \sigma_2 \! \in\!  T$, there is a unique isometry 
$f \! : \! [0,d(\sigma_1,\sigma_2)] \! \rightarrow \! T$ such
that $f(0)\!=\! \sigma_1$ and $f(d(\sigma_1,\sigma_2))\! =\! \sigma_2$. In this case, we set 
$\llb \sigma_1,\sigma_2\rrb \! :=\! f([0,d (\sigma_1,\sigma_2)])$. 

\smallskip

\noi
(b)~~For any continuous injective function 
$q: [0, 1] \! \rightarrow \! T$,  $q([0,1]) \! = \! \llb q(0), q(1)\rrb$.

\smallskip

\noi
When a point $\rho \! \in\!  T$ is distinguished, $(T, d, \rho)$ is said to be a \textit{rooted} real tree, $\rho$ being the \textit{root} of $T$. 
Among connected metric spaces, real trees are characterized by the so-called {\em four-point inequality }: 
we refer to Evans \cite{evans05} or to Dress, Moulton \& Terhalle \cite{DMT96} for a detailed account on this property. 
Let us briefly mention that the set of (pointed) isometry classes of compact rooted real trees can be equipped with the (pointed) Gromov--Hausdorff distance which makes it into a Polish space: see Evans, Pitman \& Winter \cite{EPW}, Theorem 2, for more detail on this intrinsic point of view that we do not adopt here.

\paragraph{Coding of real trees.} Real trees can be constructed through continuous functions. Recall that $e$ stands for the canonical process on $\bC(\bbR_+, \bbR_+)$. We assume here that $e$ 
has a compact support, that $e_0\! = \! 0$ and that $e$ is not identically null. Recall from (\ref{lifetimedef}) the definition of its lifetime $\zeta$. Then, our assumptions on $e$ 
entail that $\zeta \! \in \! (0, \infty)$. For $s,t\in[0,\zeta]$, we set
$$b(s, t):=\inf_{r\in[s\wedge t,s\vee t]}e_r\quad \text{ and }\quad d(s,t):=e_t +e_s-2b(s, t)\; . $$
It is easy to see that $d$ is a pseudo-distance on $[0,\zeta]$. We define the equivalence relation $\sim$ by setting 
$s\sim t$ iff $d(s,t)=0$; then we set  
\begin{equation}\label{defBT}
\cT:=[0,\zeta]/ \sim \; .
\end{equation}
The function $d$ induces a distance on the quotient set $\cT$ that we keep denoting by $d$ for simplicity. 
We denote by $p: [0,\zeta] \! \rightarrow  \! \cT$
the canonical projection. Clearly $p$ 
is continuous, which implies that $(\cT,d)$ is a compact metric space.
Moreover, it is shown that $(\cT,d)$ is a real tree (see Duquesne \& Le Gall \cite{Duquesne05}, Theorem 2.1, for a proof). 
We take $\rho=p(0)$ as the \textit{root} of $\cT$. The \textit{total height} and the \textit{diameter} of $\cT$ are thus given by  
\begin{equation}
\label{htdmdef}
\Gamma =\max_{\sigma \in \cT} d (\rho, \sigma)= \max_{t\ge 0} e_t \quad \textrm{and} \quad 
D =\max_{\sigma, \sigma^\prime \in \cT} d ( \sigma, \sigma^\prime)= \max_{s,t\ge 0} \big(e_t+e_s\! -\! 2 b(s, t) \big)  .
\end{equation}
We also define on $\cT$ a finite measure ${\rm m}$, called  the \textit{mass measure}, that is the pushforward measure of the Lebesgue measure on $[0, \zeta]$ by the canonical projection $p$. Namely, for all continuous functions $f:\cT\! \rightarrow \! \bbR_+$, 
\begin{equation}
\label{massmeadef}
\int_{\cT} \! f( \sigma) \, {\rm m} (d\sigma) =\int_0^\zeta \!\! f(p(t)) \, dt\; .
\end{equation}
Note that ${\rm m} (\cT)= \zeta$.

\paragraph{Brownian tree. } The random rooted compact real tree $(\cT,d, \rho)$ coded by $e$ under the normalized Brownian excursion measure $\bN_{\! {\rm nr}}$ defined in (\ref{Nnorm_br})
is the \textit{Brownian tree}. Here, we recall some properties of the Brownian tree. To that end, 
for any $\sigma \! \in \! \cT$, we denote by ${\rm n} (\sigma)$ the number of connected components of 
the open set $\cT \backslash \{ \sigma\}$. Note that ${\rm n} (\sigma)$ is possibly infinite. We call this number the \textit{degree} of $\sigma$. We say that $\sigma$ is a \textit{branch point} if 
${\rm n} (\sigma) \! \ge\!  3$ and that $\sigma$ is a \textit{leaf} if ${\rm n} (\sigma) \! =\!  1$. 
We denote by $\mathtt{Lf} (\cT)\! :=\! \big\{ \sigma \! \in \! \cT: {\rm n} (\sigma) \! =\! 1 \big\} $ the \textit{set of leaves} of $\cT$. Then the following holds true: 
\begin{equation}
\label{CRTprop_br}
\textrm{$\bN_{\! {\rm nr}}$-a.s.} \quad \forall \sigma \in \cT, \quad   {\rm n} (\sigma) \in \{ 1,2,3\}, 
 \quad  \textrm{${\rm m}$ is diffuse} \quad \textrm{and} \quad {\rm m} \big( \cT \backslash \mathtt{Lf} (\cT) \big) = 0 \; ,
\end{equation}
where we recall from (\ref{massmeadef}) that ${\rm m}$ stands for the mass measure. The Brownian tree has therefore only binary branch points ({\em i.e.}\,branch points of degree $3$). The fact that the mass measure is diffuse and supported by the set of leaves makes the Brownian tree a \textit{continuum random tree} according to Aldous' terminology (see Aldous \cite{aldcrt3}). For more detail on \eqref{CRTprop_br}, see for instance Duquesne \& Le Gall \cite{Duquesne05}. 

 The choice of the normalizing constant $\sqrt{2}$ for the underlying Brownian motion $X$ is motivated by the following fact: 
let $T^*_n$ be uniformly distributed on the set of rooted \textit{planar} trees with $n$ vertices. 
We view $T_n^*$ as a graph embedded in the clockwise oriented upper half-plane, whose edges are segments of unit length and whose root is at the origin. Let us consider a particle that explores $T^*_n$ as follows: it starts at the root and then it moves continuously on the tree at unit speed from the left to the right, backtracking as less as possible. During this exploration the particle visits each edge exactly twice and its journey lasts 
$2(n\! -\! 1)$ units of time. For all $t\! \in \! [0, 2(n\! -\! 1)]$, we denote by $C^{(n)}_t$ the distance between the root and the position of the particle at time $t$. The process 
$(C_{t}^{(n)})_{ t\in[0,2(n\!-\!1)]}$ is called the \textit{contour process} of $T^*_n$. 
Following an idea of Dwass \cite{Dw}, we can check that the contour process 
$(C_{t}^{(n)})_{t\in[0,2(n\!-\!1)]}$ is distributed as the (linear interpolation of the) simple random walk starting from $0$, 
conditioned to stay nonnegative on $[0,2(n\!-\!1)]$ and to hit the value $\! -1$ at time $2n\! -\! 1$. Using a variant of Donsker's 
invariance principle, the rescaled contour function $(n^{-\frac{1}{2}}C_{2(n\!-\!1)t}^{(n)})_{ t\in[0,1]}$ converges in law 
towards $e$ under $\bN_{\! {\rm nr}}$: see for instance Le Gall \cite{LGCornell}. Thus, 
$$ n^{-\frac{1}{2}}D^*_n\quad \overset{\textrm{(law)}}{-\!\!\! -\!\!\!\longrightarrow} \quad D \quad \textrm{under}\quad \bN_{\! \textrm{nr}},  $$
where $D_n^*$ stands for the diameter of $T_n^*$ and  $D$ is the diameter of the Brownian tree given by (\ref{htdmdef}). 
\begin{rem}
\label{normCRT} In the first paragraph of the introduction, we have introduced the random tree $T_n$, which is uniformly distributed on the set of rooted labelled  trees with $n$ vertices. The law of $T_n$ is therefore distinct from that of $T^*_n$, which is uniformly distributed on the set of rooted planar trees with $n$ vertices. 
Aldous \cite{aldcrt3} has proved that the tree $T_n$, whose graph distance is rescaled by a factor $n^{-\frac{1}{2}}$, converges to the tree coded by 
$\sqrt{2} e$ under $\bN_{\! \textrm{nr}}$. Thus, 
\begin{equation}
\label{identiflaw}
\Delta \quad \overset{\textrm{(law)}}{=} \quad \sqrt{2} D \quad \textrm{under} \quad \bN_{\! \textrm{nr}}\; . 
\end{equation}
See Remark \ref{DeltaD} below. \cq
\end{rem}

In this article, we prove the following result that characterizes the joint law of the height  and  diameter of the Brownian tree. 
\begin{thm}
\label{Laphtdmth} 
Recall from 
(\ref{Nnorm_br}) the definition of $\bN_{\! {\rm nr}}$ and recall from (\ref{htdmdef}) the definitions of $\Gamma$ and $D$.  
We set  
\begin{equation}\label{eq: def_bL}
\forall \lam, y, z \in (0, \infty) , \quad {\rm L}_{\lambda} (y,z):= \frac{1}{2\sqrt{\pi}} \int_0^\infty e^{-\lambda r} 
r^{-\frac{3}{2}} \, \bN_{\! {\rm nr}} 
\big(\, r^{\frac{1}{2}} D\! >\! 2y \, ;\,  r^{\frac{1}{2}} \Gamma \! >\! z \big) \, dr \; .
\end{equation}
Note that 
\begin{equation}\label{scalebL_br}
\forall \lam, y, z \in (0, \infty) , \quad {\rm L}_1(y, z)= \lambda^{-\frac{1}{2}} {\rm L}_\lambda \big( \lambda^{-\frac{1}{2}} y\, ,\,  \lambda^{-\frac{1}{2}} z \big)   \; .
\end{equation}
Then, 
\begin{equation}
\label{Lreecrit_br}
{\rm L}_1(y, z)
=  \coth ( y  \! \vee \!   z) -1 
-\frac{_1}{^4} \un_{\{ z \le 2y\}}  \frac{\sinh (2q) -2q}{\sinh^{4} (y)} \; , 
\end{equation} 
where $q\! = \! y \! \wedge \! (2y\! -\! z)$. In particular, 
this implies that 
\begin{equation}
\label{htLapl}
\forall \lambda , z \in (0, \infty), \quad  {\rm L}_\lambda (0, z)= \sqrt{\lambda} \coth (z\sqrt{\lambda}) -\sqrt{\lambda} 
\end{equation}
and 
\begin{equation}
\label{dmLapl}
\forall \lambda , y \in (0, \infty), \quad  {\rm L}_\lambda (y, 0) 
=  \sqrt{\lambda} \coth ( y \sqrt{\lambda}) -\sqrt{\lambda} -\sqrt{\lambda} \, \frac{\sinh (2y\sqrt{\lambda}) -2y\sqrt{\lambda}}{4\sinh^{4} (y\sqrt{\lambda})} \; .
\end{equation}
\end{thm}
\begin{cor}
\label{expanhtdm} For all $y,z\! \in \! (0, \infty)$, we set 
\begin{equation}
\label{defdelta}
\rho=z\vee \frac{_y}{^2} \quad \text{ and }\quad \delta  =  \big(\frac{_{2(y-z)}}{^y}\vee 0\big)\wedge 1\; .
\end{equation} 
Then we have 
\begin{align}
\label{jtlaw} 
&\bN_{\! {\rm nr}} \big(D \! >\! y \, ; \, \Gamma \! >\! z \big)=2\sum_{n\ge 1}\(2n^2\rho^2\! -\! 1\)e^{-n^2\rho^2} +\\
&\frac{1}{6}\! \sum_{n\ge 2} n(n^2\! -\! 1) 
\Big[ \big[ (n\! +\! \delta)^2y^2 \! -\! 2 \big]e^{-\frac{1}{4} (n+\delta )^2y^2 } 
\! -\! \big[ (n\! -\! \delta)^2y^2 \! -\! 2 \big]e^{-\frac{1}{4} (n-\delta)^2y^2 } + \delta y ( n^3y^3\! -\! 6ny)e^{-\frac{1}{4}n^2y^2}
 \Big] \nonumber
\end{align}
and 
\begin{align}
\label{jtlaw'}
&\bN_{\! {\rm nr}} \big(D \! \le\! y \, ; \, \Gamma \! \le\! z \big)=\frac{4\pi^{5/2}}{\rho^3}\sum_{n\ge 1}n^2e^{-n^2\pi^2/\rho^2}-\\ \notag
&\qquad\qquad\qquad\frac{32\pi^{3/2}}{3}\sum_{n\ge 1}n\sin(2\pi n\delta)\Big(\frac{2}{y^5}(2a^2_{n, y}-9a_{n, y}+6)-\frac{3\delta^2-1}{y^3}(a_{n, y}-1)\Big)e^{-a_{n, y}}+\\ \notag
&\qquad\qquad\qquad\frac{16\pi^{1/2}}{3}\sum_{n\ge 1}\delta\cos(2\pi n\delta)\Big(\frac{1}{y^3}(6a^2_{n, y}-15a_{n, y}+3)-\frac{\delta^2-1}{2y}a_{n, y}\Big)e^{-a_{n, y}}+\\ \notag
&\qquad\qquad\qquad\frac{16\pi^{1/2}}{3}\sum_{n\ge 1}\delta\Big(\frac{1}{y^3}(4a^3_{n, y}-24a^2_{n, y}+27a_{n, y}-3)+\frac{1}{2y}(2a^2_{n, y}-3a_{n, y})\Big)e^{-a_{n, y}},
\end{align}
where we have set  $a_{n,y}= 4 (\pi n/y)^2$ for all $y \! \in \! (0, \infty)$ and for all $n\! \ge \! 1$ to simplify notation. In particular, \eqref{jtlaw} implies 
\begin{equation}
\label{htBronor_bt}
\bN_{\! {\rm nr}} \big(\Gamma \! >\! y \big)=2\sum_{n\ge 1}\(2n^2y^2-1\)e^{-n^2y^2} , 
\end{equation}
and
\begin{equation}
\label{diaBronor_br}
\bN_{\! {\rm nr}} \big( D \! >\! y\big) = 
\sum_{n\ge 2} (n^2-1)\big(\frac{_1}{^6}n^4y^4-2n^2y^2+2\big)e^{-n^2y^2/4}  .
\end{equation}
On the other hand,  \eqref{jtlaw'} implies 
\begin{equation}
\label{bhtBronor_bt} 
\bN_{\! {\rm nr}} \big(\Gamma \! \le \! y \big)=
\frac{4\pi^{5/2}}{y^{3}}
\sum_{n \ge 1} n^2e^{-n^2\pi^2/y^2}, 
\end{equation}
and
\begin{equation}
\label{bdiaBronor_bt} \bN_{\! {\rm nr}} \big( D \! \le \! y\big)= 
\frac{\sqrt{\pi}}{3} \sum_{n\ge 1} \Big(  \frac{8}{y^3} \big( 24a_{n,y} -36 a_{n,y}^2 +8a_{n,y}^3 \big) +\frac{16}{y} a_{n,y}^2 \Big) e^{-a_{n,y}} \; .
\end{equation}
Thus the law of $D$ under $\bN_{\! {\rm nr}} $ has the following density: 
\begin{eqnarray}
\label{diadens}
\!\!\!\!\!\!\!\!\!\!\!\!\!\!  f_D(y)&= & \frac{1}{12} \sum_{n\ge 1}    \big(  n^8y^5-n^6y^3(20+y^2) +20 n^4y(3+y^2) -60n^2y \big) e^{-n^2y^2/4}\\
\label{bdiadens} &= & 
\frac{2\sqrt{\pi}}{3} \sum_{n\ge 1} \Big(  \frac{16}{y^4} (4a_{n,y}^4 -36a_{n,y}^3+75a_{n,y}^2 -30a_{n,y}) + \frac{8}{y^2} (2a_{n,y}^3 -5a_{n,y}^2 )\Big)e^{-a_{n,y}} \; . 
\end{eqnarray}
\end{cor}

\begin{rem}
\label{Jacobir} 
We derive \eqref{jtlaw'} from \eqref{jtlaw} using the following identity on the theta function due to Jacobi (1828), which is a consequence of Poisson summation formula:
\begin{equation}
\label{Jacobi}
\forall t\in (0, \infty), \forall x,y\in \bbC, \quad \sum_{n\in \bbZ} e^{- (x+n)^2t-2\pi i ny} = e^{2\pi ixy}\(\frac{\pi}{t}\)^{\frac{1}{2}}\sum_{n\in \bbZ} e^{-\frac{\pi^2(y+n)^2}{t}+2\pi i nx} \; .
\end{equation} 
See for instance Weil \cite{Weil}, Chapter VII, Equation (12).  Not surprisingly,  \eqref{Jacobi} can also be used to derive 
(\ref{bhtBronor_bt}) from (\ref{htBronor_bt}), to derive (\ref{bdiaBronor_bt}) from (\ref{diaBronor_br}), or to derive (\ref{bdiadens}) from (\ref{diadens}).\cq 
\end{rem}

\begin{rem}
\label{DeltaD} 
We obtain (\ref{diadens}) (resp. (\ref{bdiadens})) by differentiating (\ref{diaBronor_br}) (resp.  (\ref{bdiaBronor_bt})). 
By (\ref{identiflaw}), we have 
$$ \forall y\in (0, \infty), \quad f_\Delta (y)= \frac{1}{\sqrt{2}} f_D \Big( \frac{y}{\sqrt{2}}\Big) \; ,$$ 
which immediately entails (\ref{eq: sz}) from (\ref{bdiadens}), since $a_{n, y/\sqrt{2}}= 8(\pi n/y)^2= b_{n, y}$. \cq 
\end{rem}

\begin{rem}
Recall that $\Gamma=\max_{t\ge 0}e_t$.  Equations (\ref{htBronor_bt}) and (\ref{bhtBronor_bt}) are consistent with previous results on the distribution of the maximum of Brownian excursion: see for example Chung \cite{Chung}, though we need to keep in mind the difference between $\bN_{\! {\rm nr}}$ and $\bN_{\! {\rm Ito}}^+$, as explained in Remark \ref{nonstandard}. \cq
\end{rem}

\medskip

\noi\emph{Acknowledgements. } The author is deeply grateful to Thomas Duquesne for suggesting this problem and for those fruitful discussions. The author would also like to thank Romain Abraham and Svante Janson for a careful reading of an earlier version.

\section{Preliminaries}
\label{Prelsec}
\paragraph{A geometric property on diameters of real trees.} We begin with a simple observation on the total height and diameter of a real tree. 
\begin{lem}\label{lem: dm_br}
Let $(\rT, d, \rho)$ be a compact rooted real tree. Then $\Gamma \! \le \! D \! \le \! 2\Gamma$, where 
$$\Gamma\! =\! \sup_{u\in \rT} d(u, \rho) \quad \textrm{and} \quad  D\! = \! \sup_{u,v\in \rT} d(u,v) \; .$$ 
Moreover, there exists a pair of points $u_0, v_0 \! \in \! \rT$ with maximal distance. Namely, 
\begin{equation}\label{eq: bDT}
d(u_0, v_0)
=\sup_{u, v\in \rT}d(u, v) = D \; .
\end{equation}
Without loss of generality, we assume that $d(u_0, \rho) \! \ge \! d(v_0, \rho)$. Then the total height of $\rT$ is attained at $u_0$. Namely 
\begin{equation}\label{eq: DTGamT_br}
d( u_0, \rho)=
\sup_{u\in T}d(u, \rho)= \Gamma \; .
\end{equation}
\end{lem}

\begin{figure}[htp]
\centering
\includegraphics[scale=.9]{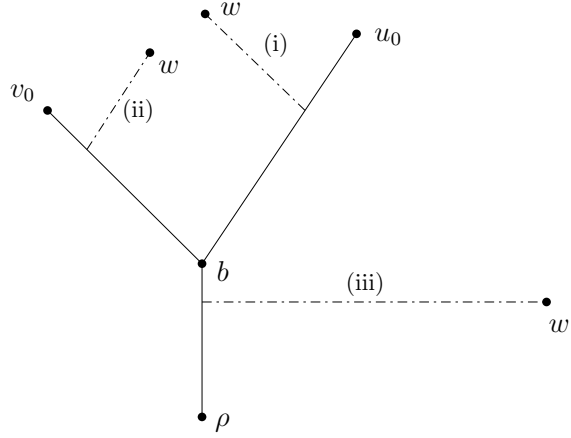}
\caption{\label{fig: 1-br}Three possibilities for $w$}
\end{figure}

\noi
\textbf{Proof. } Let $u,v\! \in \! \rT$. Recall from the definition of real trees (given in the introduction) that $\llb u, v\rrb$ stands for the unique geodesic path between $u$ and $v$. To simplify notation, we set $h(u)\! :=\! d(u, \rho)$ for $u \! \in \! \rT$.
The \textit{branch point} $u \! \wedge \! v$ of $u$ and $v$ is the unique point of $\rT$ satisfying 
$$\llb \rho, u\wedge v \rrb=\llb \rho, u\rrb\cap\llb \rho, v\rrb \; .$$ 
Then, we easily check  that
$$ d(u,v)= d(u, u\wedge v)+ d(u\wedge v, v)= h(u)+h(v)-2h(u\wedge v) \; .$$
The triangle inequality easily implies that $D \! \le \! 2 \Gamma$ while the inequality $\Gamma \! \le \! D$ is a direct consequence of the definitions. As $d: \rT^2 \rightarrow \bbR_+$ is continuous and $\rT$ is compact, there exists a pair of points $u_0, v_0 \! \in \! \rT$ such that (\ref{eq: bDT}) holds true. 
 To prove \eqref{eq: DTGamT_br}, we argue by contradiction: we assume that there exists $w \! \in \! \rT$ such that $h(w)\! >\! h(u_0)$. 
Let us write $b \! :=\! u_0 \! \wedge \! v_0$. 
Here we enumerate the three possible locations of $w$. See Figure \ref{fig: 1-br}.
\begin{enumerate}[(i)]
\item
Suppose that $w\wedge u_0\in \llb u_0, b\rrb$.
By hypothesis, we have $h(w)>h(u_0)$. In other words, 
$$
h(w)=d(w, b)+h(b)>h(u_0)=d(u_0, b)+h(b).
$$
Thus, $d(w, b)>d(u_0, b)$ and 
$$
d(w, v_0)=d(w, b)+d(b, v_0)
>d(u_0, b)+d(b, v_0)=d(u_0, v_0),
$$
which contradicts \eqref{eq: bDT}.

\item 
Suppose that $w\wedge v_0\in \llb v_0, b\rrb$.
In this case, we have
$$
h(w)=d(w, b)+h(b)>h(u_0)\ge h(v_0)=d(v_0, b)+h(b).
$$
Then $d(w, b)>d(v_0, b)$ and 
$$
d(w, u_0)=d(w, b)+d(b, u_0)
>d(v_0, b)+d(b, u_0)=d(u_0, v_0).
$$
This again contradicts \eqref{eq: bDT}.

\item
Suppose that $w\wedge u_0\in \llb \rho, b\rrb$. 
Then we deduce from 
$$
h(w)=d(w, w\wedge u_0)+h(w\wedge u_0)>h(u_0)=d(u_0, w\wedge u_0)+h(w\wedge u_0)
$$
that $d(w, w\wedge u_0)>d(u_0, w\wedge u_0)$. Note that in this case $w \wedge u_0= w\wedge v_0$.  
Therefore,
\begin{align*}
d(w, v_0)=d(w, w\wedge v_0)+d(w\wedge v_0, v_0) & = d(w, w\wedge u_0)+d(w\wedge u_0, v_0) \\ 
& >d(u_0, w\wedge u_0)+d(w\wedge u_0, v_0)\\
&> d(u_0, b)+d(b, v_0)=d(u_0, v_0),
\end{align*}
which contradicts \eqref{eq: bDT}.
\end{enumerate}
In short, there exists no $w\! \in \! \rT$ such that $h(w)=d(w, \rho) >h(u_0)= d(u_0, \rho)$, which entails \eqref{eq: DTGamT_br}.  \cqfd

\paragraph{Williams' decomposition of Brownian excursions.} Let us recall the classical result of Williams' path decomposition of Brownian excursions (see for instance 
Revuz \& Yor \cite{yorbook} Chapter XII Theorem 4.5).  
Define 
\begin{equation}
\label{taustar}
\tau_*:=\inf\{t\! >\! 0: e_t\! =\! \max_{s\ge 0} e_s\}\; .
\end{equation}
Then, 
\begin{equation}\label{unique}
\bN\text{-a.e. } ( \bN_{\! {\rm nr}}\text{-a.s.}),\, \tau_* \text{\emph{ is the unique time at which $e$ reaches its maximum value.} }
\end{equation}
Recall from (\ref{htdmdef}) the definition of the total height $\Gamma$ of the Brownian tree coded by $e$. Then, we have $\Gamma \! = \! e_{\tau_*}$.

We also recall the distribution of $\Gamma$ under $\bN$: 
\begin{equation}\label{eq: d_Gam}
\bN \big(\Gamma \! \in \! dr\big)=\frac{dr}{r^2} \; .
\end{equation}
See Revuz \& Yor \cite{yorbook} Chapter XII Theorem 4.5 combined with Remark \ref{nonstandard}.  

Williams's decomposition entails that there is a regular version of the family of conditioned laws $\bN (\, \cdot \, | \, \Gamma \! = \! r)$, $r\! >\! 0$. Namely, $\bN (\, \cdot \, | \, \Gamma \! = \! r)$-a.s.~$\Gamma \! = \! r$,
$r\mapsto \bN (\, \cdot \, | \, \Gamma \! = \! r)$ is weakly continuous on $\bC(\bbR_+, \bbR_+)$ and 
\begin{equation}
 \bN = \int_0^\infty  \!\!\! \bN(\Gamma \! \in \! dr)\,  \bN (\, \cdot \, | \, \Gamma \! = \! r) \; .
 \end{equation}

Let $Z\! =\! (Z_t)_{ t\ge 0}$ be a continuous process defined on the probability space $(\Omega, \cF, \bP)$ such that $\tfrac{1}{\sqrt 2}Z$ is distributed as a Bessel process of dimension $3$ starting from $0$. Let $\tau_r=\inf\{t\! >\! 0: Z_t \! =\! r\}$ be the hitting time of $Z$ at level $r\! \in \! (0, \infty)$. 
We recall that 
\begin{equation}
\label{hitBess}
\forall \lambda \in \bbR_+, \quad \bE \big[ e^{-\lambda \tau_r}\big] = \frac{r\sqrt{\lam}}{\sinh(r\sqrt{\lam})} \; .
\end{equation}
See Borodin \& Salminen \cite{BoroSalm} Part II, Chapter 5, Section 2, Formula 2.0.1, p.~463, where we let $x$ tend to $0$ and take $\alpha\! = \! \lambda$ and 
$z\! = \! r/\sqrt{2}$, since $Z\! =\!  \sqrt{2}R^{(3)}$.  

We next introduce the following notation
$$
\overleftarrow{e}(t)=e_{(\tau_*-t)_+}; \quad \overrightarrow{e}(t)=e_{\tau_* +t}, \quad t\ge 0.
$$
where $(\cdot)_+$ stands for the positive part function. \textit{Williams' decomposition} of Brownian excursion asserts that 
\begin{itemize}
\item[] \textit{for all $r\! \in \! (0, \infty)$, under $\bN(\, \cdot \, | \, \Gamma\! =\! r)$, the two processes $\overleftarrow e$ and $\overrightarrow e$ are distributed as two independent copies of $(Z_{(\tau_r-t)_+})_{ t\ge 0}$. }
\end{itemize}
As a combined consequence of this decomposition and  (\ref{hitBess}), we have
\begin{equation}
\label{eq: lap_tau}
\forall r\! \in \! (0, \infty), \quad \bN\big(e^{-\lam\zeta}|\, \Gamma=r\big)=\bE \big[e^{-\lam \tau_r} \big]^2=\(\frac{r\sqrt{\lam}}{\sinh(r\sqrt{\lam})}\)^2, 
\end{equation}
where we recall that $\zeta$ stands for the lifetime of the excursion. Therefore,
$$
\bN\(e^{-\lam\zeta}\indi_{\{\Gamma>a\}}\)=\int_a^\infty \!\!\! \bN\big(e^{-\lam \zeta}|\, \Gamma=r \big)\bN(\Gamma\in dr)
=\int_a^\infty \!\!\!\frac{\lam dr}{\sinh^2(r\sqrt{\lam})}=\sqrt{\lam}\coth(a\sqrt{\lam})-\sqrt{\lam}\, ,
$$
by \eqref{eq: d_Gam} and \eqref{eq: lap_tau}. Combined with the fact that $\bN(1-e^{-\lam\zeta})=\sqrt{\lam}$, this entails that
\begin{equation}
\label{w_br}
\bN\Big(1-e^{-\lam \zeta}\indi_{\{\Gamma\le a\}}\Big)=\sqrt{\lam}\coth(a\sqrt{\lam}).
\end{equation}
This equation is used in the proof of Theorem \ref{Laphtdmth}.

\paragraph{Spinal decomposition}
Let us interpret Williams' decomposition in terms of a Poisson decomposition of the Brownian excursion. To that end, we need the following notation.  
Let $h \! \in \! \bC(\bbR_+, \bbR_+)$ have compact support. We assume that $h(0) \! >\! 0$.  For all $s\! \in \! \bbR_+$, we set 
$\underline{h}(s) \! =\! \inf_{0\le u\le s}h(u)$. Let  $(l_i, r_i)$, $i\! \in \! \cI(h)$ be the excursion intervals of $h\! -\! \underline{h}$ away from $0$; namely, 
they are the connected components of the open set $\{s\ge 0: h(s)\! -\! \underline{h}(s)\! >\! 0\}$. For all $i \! \in \! \cI(h)$, we next set 
$$
h^i(s)=\big(h-\underline{h}\big)\big((l_i+s)\wedge r_i\big), \quad s\ge 0 , 
$$
which is the excursion of $h\! -\! \underline{h}$ corresponding to the interval $(l_i, r_i)$. Then we set 
$$
\cP(h)=\sum_{i\in \cI(h)}\delta_{(h(0)-h(l_i), \,h^i)}, 
$$
that is a point measure on $\bbR_+ \! \times \! \bC(\bbR_+, \bbR_+) $. We define 
\begin{equation}\label{eq: spinal_e}
\cQ:=\cP(\overrightarrow{e})+\cP(\overleftarrow{e})=:\sum_{j\in \cJ}\delta_{(s_j, \,e^j)} \; .
\end{equation}
We also introduce  
for all $t \! \in \! (0, \infty)$ the following notation
\begin{equation}
\label{NotaNt}
\bN_t = \bN \big( \cdot \, \cap \, \{ \Gamma \! \le \! t \} \big) \; .
\end{equation}
The following lemma is the special case of a general result due to Abraham \& Delmas \cite{AbDe09}. 
\begin{lem}[Proposition 1.1, Abraham \& Delmas \cite{AbDe09}]
\label{lem: williams_br}
Let $r\! \in \! (0, \infty)$. Then, $\cQ$ under $\bN(\, \cdot \, |\, \Gamma \! =\! r)$ is a Poisson point measure on $\bbR_+ \! \times \! \bC(\bbR_+, \bbR_+) $
with intensity measure $2\!\cdot\!\indi_{[0, r]}(t) dt\, \bN_t$.
\end{lem}

\paragraph{Interpretation in terms of the Brownian tree and consequences.} 
Let us interpret $\cQ$ in terms of the Brownian tree $\cT$ coded by the Brownian excursion $e$. 
Recall that $p: [0, \zeta] \! \to \! \cT$ stands for the canonical projection and recall that $\rho\! =\! p(0)$ is the root of $\cT$. 
The point $p(\tau_*)$ is the (unique) point of $\cT$ that attains the total height: $d(\rho, p(\tau_*)) \! = \! \Gamma$. 

Denote by $\cT_{j'}^o$, $j'\! \in \! \cJ'$, the connected components of $\cT \backslash \llb \rho, p(\tau_*)\rrb$. For all 
$j'\! \in \! \cJ'$, there exists a unique point $\sigma_{j'} \! \in \! \llb \rho, p(\tau_*)\rrb$ such that $\cT_{j'}\! :=\! \cT_{j'}^o\cup \{ \sigma_{j'}\}$ is the closure of $\cT^o_{j'}$ in $\cT$. Recall the notation $\cJ$ from (\ref{eq: spinal_e}). It is not difficult to see that $\cJ' $ is in one-to-one correspondence with $\cJ$. Moreover, after a re-indexing, we can suppose that $d(p(\tau_*), \sigma_j)\! = \! s_j$ and that $(\cT_j, d, \sigma_j)$ is the real tree coded by the excursion $e^j$, for each $j\in \cJ$. Then we set 
\begin{equation}
\label{defGamj}
\forall j\in \cJ, \quad \Gamma_j:= \max_{s\ge 0} e^j (s) = \max_{\gamma \in \cT_j} d(\sigma_j, \gamma) \; , 
\end{equation}
that is the total height of the rooted real tree $(\cT_j, d, \sigma_j)$. We claim that 
\begin{equation}\label{eq: bdcj}
\bN\text{-a.e.} \quad D=\sup_{j\in \cJ} \, (s_j+\Gamma_j) \; .
\end{equation}
\noi
\textit{Proof of (\ref{eq: bdcj}).} First observe that for all $t \in \! (0, \infty)$, $\bN_t$ is an infinite measure because $\bN$ is infinite and because 
$\bN (\Gamma \! >\! t)\! = \! 1/t$ by (\ref{eq: d_Gam}). By Lemma \ref{lem: williams_br}, $\bN$-a.e.~the closure of the 
set $\{ s_j\, ; \, j\! \in \! \cJ\}$ is $[0, \Gamma]$. This entails that 
\begin{equation}
\label{stepun}
\bN\text{-a.e.} \quad \Gamma =\sup_{j\in \cJ}s_j\le \sup_{j\in \cJ} \, (s_j+ \Gamma_j) \; .
\end{equation}  
Next, for all $j\! \in \! \cJ$, there exists $\gamma_j \! \in \! \cT_j$ such that $d(\sigma_j, \gamma_j)\! = \! \Gamma_j$. Then observe that 
\begin{equation}
\label{stepdeu}
d(p(\tau_*), \gamma_j)= d(p(\tau_*), \sigma_j)+ d(\sigma_j, \gamma_j)= s_j+ \Gamma_j \; .
\end{equation}
Note that Lemma \ref{lem: dm_br} and \eqref{unique} imply that $D\! = \! \max_{\gamma \in \cT} d(p(\tau_*), \gamma)$. Comparing this with (\ref{stepdeu}), we get 
\begin{equation}
\label{steptroi}
D \ge \sup_{j\in \cJ} \, (s_j+ \Gamma_j) \; .
\end{equation}
On the other hand, there exists $\gamma^*\! \in \! \cT$ such that $D\! = \! \max_{\gamma \in \cT} d(p(\tau_*), \gamma)\! = \! d(p(\tau_*), \gamma^*)$ by Lemma \ref{lem: dm_br}. If $\gamma^* \! \notin \! \llb \rho, p(\tau_*)\rrb$, then there exists $j^*\! \in \! \cJ$ such that $\gamma^*\! \in \! \cT_{j^*}$. In consequence, we have $D= d(p(\tau_*), \gamma^*) \le s_{j^*}+ \Gamma_{j^*}$, and then 
$D \! =\! \sup_{j\in \cJ} \, (s_j+ \Gamma_j)$ when compared with \eqref{steptroi}. If $\gamma^* \! \in \! \llb \rho, p(\tau_*)\rrb$, then (\ref{stepun}) implies that $\gamma^*\! = \! \rho$ and $D\! = \! \Gamma$. In both cases  (\ref{eq: bdcj}) holds true. \cqfd

\medskip

We next denote by $\zeta_j$ the lifetime of $e^j$ for all $j\! \in \cJ$ and  prove the following statement.  
\begin{equation}
\label{lifetimej}
\bN\text{-a.e.} \quad \sum_{j\in \cJ}\zeta_j=\zeta \; .
\end{equation} 
\noi
\textit{Proof of (\ref{lifetimej}).} Let   
$\sigma \! \in \! \llb p(\tau_*) , \rho \rrb $ be distinct from $p(\tau_*)$ and $\rho$. Then ${\rm n} (\sigma)\!  \ge \! 2$ and 
$\sigma$ is not a leaf of $\cT$. Recall from (\ref{massmeadef}) the definition of the mass measure ${\rm m}$ and 
recall from (\ref{CRTprop_br}) that $\bN_{\! {\rm nr}}$-a.s.~${\rm m}$ 
is diffuse and supported on the set of leaves of $\cT$. By (\ref{zetades_br}), this property also holds true $\bN$-almost everywhere and we thus get 
\begin{equation*} 
\bN\text{-a.e.} \quad {\rm m} \big( \llb p(\tau_*) , \rho \rrb \big)= 0\; .
\end{equation*}
Recall that $\cT_j^o$, $j\! \in \! \cJ$, are the connected components of 
$\cT \backslash \llb \rho, p(\tau_*)\rrb$. Thus, 
\begin{equation} 
\label{kkugzdec}
\bN\text{-a.e.} \quad {\rm m} (\cT) = {\rm m}\big( \llb p(\tau_*) , \rho \rrb \big)+ \sum_{j\in \cJ} {\rm m} \big( \cT_j^o\big)=\sum_{j\in \cJ} {\rm m} \big( \cT_j^o\big).
\end{equation}
Recall that $\cT_j= \cT_j^o \cup \{ \sigma_j\}$ and that ${\rm m}$ is $\bN$-a.e.~diffuse, which entails ${\rm m} (\cT_j)= {\rm m} (\cT_j^o)$, for all $j\! \in \! \cJ$. Moreover, since $(\cT_j, d, \sigma_j)$ is coded by the excursion $e^j$, we have $\zeta_j= {\rm m} (\cT_j)$. For a similar reason, we also have $\zeta= {\rm m} (\cT)$. This, combined with (\ref{kkugzdec}), entails (\ref{lifetimej}). \cqfd

\section{Proof of Theorem \ref{Laphtdmth}}
\label{pfsecth}
First we note that by \eqref{echtscal_br},
\begin{equation}\label{eqq2}
{\rm L}_\lam(y, z)=\frac{1}{2\sqrt{\pi}}\int_0^\infty \!\!\! dr e^{-\lam r} r^{-\frac{3}{2}}\bN_{\! {\rm nr}} \big(r^{\frac{1}{2}}D\! >\! 2y \, ;\,  r^{\frac{1}{2}}\Gamma \! >\! z \big)
=\bN \( \! e^{-\lam \zeta}\boldsymbol{1}_{\{D> 2y; \Gamma>z\}}\).
\end{equation}
Observe that 
the scaling property \eqref{scalebL_br} is a direct consequence of the scaling property of $\bN$ (see \eqref{HscalingN_br}). 

We next 
compute the right hand side of (\ref{eqq2}). 
To that end, recall from \eqref{eq: spinal_e} the spinal decomposition of the excursion $e$ 
and recall from (\ref{defGamj}) the notation $\Gamma_j \! =\! \max_{s\ge 0} e^j (s) $, for all $j\! \in \! \cJ$; also recall that 
$\zeta_j$ stands for the lifetime of $e^j$.  
Let $r,y\! \in \! (0, \infty)$ be such that $y\! \le \! r\! \le \! 2y$. We apply successively \eqref{eq: bdcj}, (\ref{lifetimej}), Lemma \ref{lem: williams_br} and Campbell's formula for Poisson point measures and find that 
\begin{align}
\label{triton}
\bN\Big(e^{-\lam\zeta}\indi_{\{D\le 2y\}}\Big| \, \Gamma=r\Big) &= 
\bN\Big( \prod_{j\in \cJ}e^{-\lam \zeta_j}\indi_{\{s_j+\Gamma_j\le 2y\}}\Big| \, \Gamma=r\Big)  \nonumber \\
& = \exp\(\!-2\!\! \int_0^r \!\!\! dt \, \bN_t\big(1\! -\! e^{-\lam \zeta}\indi_{\{\Gamma\le 2y-t\}}\big)\)\; .
\end{align}
Recall from (\ref{NotaNt}) that $\bN_t = \bN \big( \cdot \, \cap \, \{ \Gamma \! \le \! t \} \big)$ and observe that 
\begin{equation}
\label{gloup}
\int_0^r \!\!\!\! dt \, \bN_t\big(1\! -\! e^{-\lam \zeta}\indi_{\{\Gamma\le 2y-t\}}\big)
=\int_0^{y} \!\!\!\! dt \, \bN\big((1\! -\! e^{-\lam \zeta})\indi_{\{\Gamma\le t\}}\big)  +\! \int_{y}^r \!\!\!\! dt \, \bN\big(\indi_{\{\Gamma\le t\}}\! -\! e^{-\lam \zeta}\indi_{\{\Gamma<2y-t\}}\big)\; .
\end{equation}
By \eqref{w_br} and by \eqref{eq: d_Gam}, 
\begin{equation}
\label{beurk_br}
\bN\big((1\! -\! e^{-\lam \zeta})\indi_{\{\Gamma\le t\}}\big) = \bN\big(1\! -\! e^{-\lam \zeta}\indi_{\{\Gamma\le t\}}\big) -\bN\big( \Gamma \! >\! t \big)  = \sqrt{\lambda} \coth \big( t\sqrt{\lambda} \big) -\frac{_1}{^t}\; 
\end{equation}
and 
\begin{equation}
\label{rebeurk}
\bN\big(\indi_{\{\Gamma\le t\}}\! -\! e^{-\lam \zeta}\indi_{\{\Gamma<2y-t\}}\big)\! = \! \bN\big(1\! -\! e^{-\lam \zeta}\indi_{\{ \Gamma\le 2y-t\}}\big) \!- \! \bN ( \Gamma \! >\! t ) \!  =\!  \sqrt{\lambda} \coth \big( (2y\! -\! t)\sqrt{\lambda} \big) \! -\! \frac{_1}{^t} .
\end{equation} 
Then observe that for all $\varepsilon, a \! \in (0, \infty) $ such that $\varepsilon \! < \! a$,
$$ \int_\varepsilon^a \!\!\! \big( \sqrt{\lambda} \coth(t\sqrt{\lambda}) -\frac{_1}{^t} \big) dt = \log \frac{\sinh a\sqrt{\lambda}}{a}-  
\log \frac{\sinh \varepsilon\sqrt{\lambda}}{\varepsilon} .$$
Thus, as $\varepsilon \rightarrow 0$, we get 
\begin{equation}
\label{rerebeurk}
\forall a  \in \bbR_+, \qquad \int_0^a \!\!\! \big( \sqrt{\lambda} \coth(t\sqrt{\lambda}) -\frac{_1}{^t} \big) dt =  \log \frac{\sinh a\sqrt{\lambda}}{a\sqrt{\lambda}} \; .
\end{equation}
An easy computation based on (\ref{rerebeurk}), combined with (\ref{gloup}), (\ref{beurk_br}), (\ref{rebeurk}) and (\ref{triton}), entails  
$$
\bN\Big(e^{-\lam\zeta}\indi_{\{D\le 2y\}}\Big| \, \Gamma=r\Big)=\frac{\Big(r\sqrt{\lam}\sinh\big((2y-r)\sqrt{\lam}\big)\Big)^2}{\sinh^4(y\sqrt{\lam})}.
$$
Combining this with \eqref{eq: lap_tau}, we get
\begin{equation}
\label{case1}
\forall r, y \! \in \! (0, \infty) : y\! \le \! r\! \le \! 2y, \;  \bN\big(e^{-\lam\zeta}\indi_{\{D> 2y\}}\big| \Gamma \! =\! r\big) \! =\! \(\frac{r\sqrt\lambda}{\sinh(r\sqrt\lambda)}\)^2\! -\frac{\Big(r\sqrt{\lam}\sinh\big((2y\! -\! r)\sqrt{\lam}\big)\Big)^2}{\sinh^4(y\sqrt{\lam})}.
\end{equation}
Next, let $r,y\! \in \! (0, \infty)$ be such that $r\! > \! 2y$. By Lemma \ref{lem: dm_br}, $\Gamma \! \le \! D \! \le \! 2\Gamma$. Therefore,
\begin{equation}
\label{case2}
\forall r, y \! \in \! (0, \infty) : r\! >\! 2y, \;  \bN\Big( e^{-\lam\zeta}\indi_{\{D> 2y\}}\Big|\, \Gamma=r\Big)=\bN\Big(e^{-\lam\zeta}\,\Big|\, \Gamma=r\Big)=\(\frac{r\sqrt\lambda}{\sinh(r\sqrt\lambda)}\)^2.
\end{equation}
Finally, let $r\!<\! y$. Then $\bN(e^{-\lam\zeta}\indi_{\{D> 2y\}}| \Gamma\!=\!r)\! =\! 0$, since $\Gamma \! \le \! D \! \le \! 2\Gamma$. Combining this with \eqref{case1} and \eqref{case2}, we easily obtain that
\begin{align*}
\bN\(e^{-\lambda \zeta}\boldsymbol{1}_{\{D> 2y, \Gamma>z\}}\)& =\int_{z}^\infty \!\!\! \bN\Big(e^{-\lambda\zeta}\indi_{\{D>2y\}}\Big|\, \Gamma \! =\! r\Big)\bN(\Gamma \! \in \! dr)\\
&=\int_{z\vee y}^{2y\vee z} \!\!  \bN\Big(e^{-\lambda\zeta}\indi_{\{D>2y\}}\Big|\, \Gamma \! =\! r\Big)\bN(\Gamma \! \in \! dr)+\!\! \int_{2y\vee z}^\infty \!\!
\bN\Big(e^{-\lambda\zeta}\Big|\, \Gamma \! =\! r\Big)\bN(\Gamma \! \in\!  dr)\\
&=\sqrt{\lambda}\big(\coth\big((z\! \vee \! y)\sqrt{\lambda}\big)\! -\! 1\big) \! -\! \indi_{\{z\le 2y\}}
\frac{\sqrt{\lambda}\sinh (2q\sqrt{\lambda}) \! -\! 2\lambda q}{4\sinh^4(y\sqrt{\lambda})} \; , 
\end{align*}
where we recall the notation $q=y\! \wedge \! (2y\! -\! z)$.
By \eqref{eqq2}, this concludes the proof of Theorem \ref{Laphtdmth}.

\section{Proof of Corollary  \ref{expanhtdm} }
\label{pfseccor}
We introduce the following notation for the Laplace transform on $\bbR_+$:  
for all Lebesgue integrable functions $f: \bbR_+ \to \bbR$, we set 
\begin{equation*}  
\forall \lambda \in \bbR_+, \qquad \cL_\lambda (f):= \int_0^\infty \!\!\! dx \, e^{-\lambda x}f(x), 
\end{equation*}
which is well-defined. Note that if $f, g$ are two continuous and integrable functions such that $\cL_\lambda(f)=\cL_\lambda(g)$ for all $\lambda\in [0, \infty)$, then we have 
$f\! =\! g$, by the injectivity of the Laplace transform and standard arguments. 

\medskip

For all $a,x\! \in \! (0, \infty)$, we set $f_a(x)= \frac{a}{2\sqrt{\pi}}x^{-3/2} e^{-a^2/4x}$. It is well-known that $\cL_\lambda(f_a)= 
e^{-a\sqrt{\lambda}}$ for all $\lambda \! \in \! \bbR_+$ (see for instance Borodin \& Salminen \cite{BoroSalm} Appendix 3, Particular formul{\ae}  2, p.~650). 
Then we set 
$$ g_a (x)= \partial_x f_a (x)= \frac{1}{8\sqrt{\pi}} x^{-\frac{7}{2}} e^{-\frac{a^2}{4x}} (a^3-6ax) \quad \textrm{and} \quad h_a(x) = -\partial_a f_a (x)=  
\frac{1}{4\sqrt{\pi}} x^{-\frac{5}{2}} e^{-\frac{a^2}{4x}} (a^2-2x) \; .$$
Consequently, for all $\lambda \! \in \! \bbR_+$,
\begin{equation}
\label{Lalalapl}
\cL_\lambda( g_a) = 
\lambda e^{-a\sqrt{\lambda}} \quad \textrm{and} \quad  \cL_\lambda( h_a) = 
\sqrt{\lambda}e^{-a\sqrt{\lambda}} \; .
\end{equation}
(See also Borodin \& Salminen \cite{BoroSalm} Appendix 3, Particular formul{\ae}  3 and 4, p.~650.) 
Moreover, we have the following easy bounds: for all $\lambda\in \bbR_+$,
\begin{align}
\label{bound_g}
&\cL_\lambda(|g_a|)\!\le\! \frac{1}{8\sqrt{\pi}}\int_0^\infty \!\!\! dx\, e^{-\lambda x} x^{-\frac{7}{2}} e^{-\frac{a^2}{4x}} (a^3\!+\!6ax) =  \lambda e^{-a\sqrt{\lambda}}\!+\!\frac{_6}{^a}\sqrt{\lambda}e^{-a\sqrt\lambda}\!+\!\frac{_6}{^{a^2}}e^{-a\sqrt\lambda},\\ 
\label{bound_h}
&\cL_\lambda(|h_a|)\!\le\! \frac{1}{4\sqrt{\pi}}\int_0^\infty \!\!\! dx\, e^{-\lambda x} x^{-\frac{5}{2}} e^{-\frac{a^2}{4x}} (a^2\!+\!2x)= \sqrt{\lambda}e^{-a\sqrt\lambda}\!+\!\frac{_2}{^a}e^{-a\sqrt\lambda}.
\end{align}
Let $y,z \! \in \! (0, \infty)$. Recall from \eqref{defdelta} the notation $\rho$ and $\delta$. Next set 
$$\forall n \in \bbN, \quad  u_n \! =\! \frac{_1}{^6}(n+3)(n+2)(n+1), $$ 
so that $(1\! -\! x)^{-4}\! = \! \sum_{n\ge 0} u_n x^n$, for all $x\! \in \! [0, 1)$. Then (\ref{Lreecrit_br}) implies that 
\begin{eqnarray*}
\label{expanun}
{\rm L}_1( \frac{_{_1}}{^{^2}}y , z) &= & \coth \rho \! -\! 1\! -\!  \frac{\sinh (\delta y) -\delta y }{4\sinh^4 (y/2)}= \frac{2e^{-2\rho}}{1\! -\! e^{-2\rho}}+ \frac{2e^{-2y} (e^{-\delta y}\! -\! e^{\delta y})}{(1\! -\! e^{-y})^4} + \frac{4\delta y e^{-2y}}{(1\!-\! e^{-y})^4}  \nonumber \\
& =& \sum_{n\ge 1} 2e^{-2n\rho} + \sum_{n\ge 0} 2u_n \big( e^{-(n+2+\delta)y}\! -e^{-(n+2-\delta)y} +2 \delta ye^{-(n+2)y}\big) \nonumber \\
& =& \sum_{n\ge 1} 2e^{-2n\rho} + \sum_{n\ge 2} 2u_{n-2} \big( e^{-( n+\delta)y}\! -e^{-(n-\delta)y} +2 \delta ye^{-ny}\big) \; . 
\nonumber
\end{eqnarray*}
Thus, by (\ref{scalebL_br}), we obtain that
\begin{align}
\label{expantruc}
\frac{1}{2\sqrt{\pi}} \int_0^\infty \!\!\! e^{-\lambda r} &
r^{-\frac{3}{2}} \, \bN_{\! {\rm nr}} 
\big(\, r^{\frac{1}{2}} D\! >\! y \, ;\,   r^{\frac{1}{2}} \Gamma \! >\! z \big) \, dr = {\rm L}_{\lambda} (\frac{_{_1}}{^{^2}}y,z)=  \sqrt{\lambda} {\rm L}_1( \frac{_{_1}}{^{^2}}y\sqrt{\lambda} , z\sqrt{\lambda}) 
\nonumber \\ 
= & \sum_{n\ge 1} 2\sqrt{\lambda} e^{-2n\rho\sqrt{\lambda} } + \sum_{n\ge 2} 2u_{n-2} \big( \sqrt{\lambda} e^{-( n+\delta)y\sqrt{\lambda} }\! -
\sqrt{\lambda} e^{-(n-\delta)y\sqrt{\lambda} } +2 \delta y\lambda e^{-ny\sqrt{\lambda} }\big) \nonumber \\
= & \sum_{n\ge 1} 2\cL_\lambda (h_{2n\rho}) + \sum_{n\ge 2} 2u_{n-2} \cL_\lambda \big( h_{(n+\delta)y} -h_{(n- \delta)y} +2\delta y g_{ny} \big)\; .
\end{align}
Observe that for all $r\! \in \! \bbR_+$, 
\begin{equation} 
\label{controlun}
 \sum_{n\ge 1} 2 \sup_{[0, r]}|h_{2n\rho}| + \sum_{n\ge 2} 2u_{n-2} \big( \sup_{[0, r]}|h_{(n+\delta)y}| +\sup_{[0, r]}|h_{(n - \delta)y}| +2\delta y \sup_{[0, r]}|g_{ny}| \big) < \infty \; .
\end{equation}
Then, for any $r\! \in \! \bbR_+$, we set
$$
\phi_{y, z}(r):=2\sum_{n=1}^\infty h_{2n\rho}(r)e^{-r}+\sum_{n=2}^\infty 2u_{n-2} \big( h_{(n+\delta)y}(r)e^{-r}-h_{(n-\delta)y}(r)e^{-r}+2\delta y g_{ny}(r)e^{-r} \big)\; , $$
which is well-defined and continuous thanks to (\ref{controlun}). The bounds (\ref{bound_g}) and  (\ref{bound_h})  imply that $\phi_{y,z}$ is Lebesgue integrable. 
Moreover, (\ref{expantruc}) asserts that ${\rm L}_{\lambda+1} (\frac{1}{2} y, z)\! = \! \cL_\lambda (\phi_{y,z})$. By the injectivity of the Laplace transform for continuous integrable
functions (as mentioned above), we get 
$$ \forall r \in \bbR_+, \quad \phi_{y,z} (r) = \frac{1}{2\sqrt{\pi}} e^{-r } 
r^{-\frac{3}{2}} \, \bN_{\! {\rm nr}} 
\big(\, r^{\frac{1}{2}} D\! >\! y \, ;\,   r^{\frac{1}{2}} \Gamma \! >\! z \big)\; , $$
which entails (\ref{jtlaw}) by taking $r\!=\!1$. 

Since $\Gamma \! \le \! D \! \le \! 2\Gamma$, if $z\! = \! y$, then $\bN_{\! {\rm nr}} (D\! >\! y ;\Gamma\! >\! y) \! =\!  \bN_{\! {\rm nr}} (\Gamma\! >\! y)$ and (\ref{jtlaw}) immediately implies (\ref{htBronor_bt}) because in this case $\rho=y$ and $\delta\! = \! 0$. If $z\! = \! y/2$, then $\bN_{\! {\rm nr}} (D\! >\! y ; \Gamma\! >\! y/2) \! =\!  
\bN_{\! {\rm nr}} (D\! >\! y)$, $\rho\!=\!y/2$, $\delta\! = \! 1$ and (\ref{jtlaw}) implies 
\begin{align}
\bN_{\! {\rm nr}} (D\! >\! y)= \sum_{n\ge 1} \(n^2y^2\! -\! 2\)e^{-\frac{1}{4}n^2y^2} + \frac{1}{6}\! &\sum_{n\ge 2} n(n^2\! -\! 1) \Big[ \big[ (n+1) ^2y^2 \! -\! 2 \big]e^{-\frac{1}{4} (n+1)^2y^2 }- \nonumber \\
& \! \big[ (n\! -\! 1)^2y^2 \! -\! 2 \big]e^{-\frac{1}{4} (n-1)^2y^2 } + y( n^3y^3\! -\! 6ny)e^{-\frac{1}{4}n^2y^2}
 \Big], \nonumber
\end{align} 
which entails (\ref{diaBronor_br}) by re-indexing the sums according to $e^{-n^2y^2/4}$: we leave the detail to the reader. We next derive (\ref{diadens}) by differentiating (\ref{diaBronor_br}). As mentioned in Remark \ref{Jacobir}, we use Jacobi identity (\ref{Jacobi}) to derive \eqref{jtlaw'} from \eqref{jtlaw}. The computations are long but straightforward: we leave them to the reader.
Finally, for the same reason as before, \eqref{jtlaw'} entails \eqref{bhtBronor_bt} by taking $\rho\!=\!y$ and $\delta\!=\!0$. It also entails \eqref{bdiaBronor_bt} by taking $\rho\!=\!\frac{_1}{^2}y$ and $\delta\! = \! 1$. 
Differentiating \eqref{bdiaBronor_bt} gives \eqref{bdiadens}.
This completes the proof of Corollary \ref{expanhtdm}.

{\small
\setlength{\bibsep}{.3em}

}

\end{document}